\documentclass[runningheads]{llncs}
\usepackage{amsmath}
\usepackage{booktabs} 
\usepackage{caption} 
\usepackage{subcaption} 
\usepackage{graphicx}
\usepackage{pgfplots}
\usepackage[all]{nowidow}
\usepackage[utf8]{inputenc}
\usepackage{tikz}
\usepackage{breqn}
\usetikzlibrary{er,positioning,bayesnet}
\usepackage{multicol}
\usepackage{algpseudocode,algorithm,algorithmicx}

\usepackage{enumitem}
\definecolor{blue}{HTML}{1F77B4}
\definecolor{orange}{HTML}{FF7F0E}
\definecolor{green}{HTML}{2CA02C}

\pgfplotsset{compat=1.14}

\begin{document}
\title{\line(1,0){250}\\Iterative Domain Optimization\\\line(1,0){250}}
%
%
\author{Lefgoum Raian Noufel}
%
%
\institute{Higher National School of Computer Science(ESI)
\\
BP 68M 16309 Oued Smar 
\\ Algiers, Algeria \\
\email{fl\_rayan@esi.dz}}

\maketitle              
\begin{abstract}
In this paper we study a new approach in optimization that aims to search a large domain D where a given function takes large, small or specific values via an iterative optimization algorithm based on the gradient. We show that the objective function used is not directly optimizable, however, we use a trick to approximate this objective by another one at each iteration to optimize it.
Then we explore a use case of this algorithm in machine learning
to find domains where the models output large and small values with respect of some constraints.
Experiments demonstrate the efficiency of this algorithm on five cases with models trained on the titanic dataset.

\keywords{Optimization; Gradient-based Optimization; Domain Optimization; Machine Learning.}
\end{abstract}
\section{Introduction}

Optimization is a subfield of mathematics where the goal is to find the best solution from a space of possibles solutions while continuous optimization is a subfield of optimization where the variables used in the objective function take continuous values, those objectives functions under the assumption of differentiability are easily optimizable with the gradient-based algorithm which is heavily used in many fields of science and is the key success behind deep learning models, an easy algorithm to implement since it only requires the first order derivative of the objective with respect of its variables.
Gradient-based algorithm caught the attention of many researchers which allowed the development of many variants, In the field of Deep learning, Stochastic Gradient descent performs parameters update for each training example, Mini-Batch gradient descent performs parameters update after computing the gradients of a batch-size of training example, Momentum speed up the convergence by accumulating the precedent gradient to cancel irrelevants directions, RMSProp restricts the oscillations with an adaptative learning rate for each parameter, Adaptive Moment Estimation combines both the benefits of momentum and RMSProp.

In this paper we study the iterative domain optimization, a generalization of continuous optimization to domains which aims to find a domain where a given function takes large, small or specific values. This algorithm works by approximating at each iteration of the optimization process our objective function by another differentiable function to compute its gradient. In our knowledge such algorithm haven't yet been proposed. We experimented this algorithm on complex functions like machine learning models trained on the titanic dataset. 
      
\section{Iterative Domain Optimization}

The goal of this algorithm is to find a domain D where our function takes large, small or specific values when we sample points from this domain.
\newline Let's first focus on large values, concretely the more a function f takes large values on a domain D the more its integral will be high on this domain, so a first approach would be to find a domain D which maximize the following expression: 
\begin{equation}
    \int_D f(x_1,..,x_n)
\end{equation}
\\
This expression yet is not complete, for example if our function always takes positives values then the solution of our optimization problem would be D* = ${\rm I\!R}^n$.
Instead of that we will focus on searching a domain D where our function f get a large mean ie: 
\begin{equation}
    \frac{1}{Vol(D)} \times \int_D f(x_1,..,x_n)
\end{equation}
\\
Small values are found by maximizing the following expression:
\\
\begin{equation}
    \frac{1}{Vol(D)} \times \int_D -f(x_1,..,x_n)
\end{equation}
\\
And finally if M is the average value that we want on a domain D then the expression to maximize becomes: 

\begin{equation}
-\frac{1}{Vol(D)} \left(\int_D f(x_1,..,x_n)  \; - \; M \right)^2
\end{equation}
\\ 
For what follows we will focus on searching large values. Small and specific values could be found by considering the precedents expressions, also, all the tuning parameters are supposed in ${\rm I\!R}^+$.

\subsection{The objective function} \label{obj}

Let's now define our domain as $[c_1-\sigma_1,c_1+\sigma_1] \times .. \times[c_n-\sigma_n,c_n+\sigma_n]$ where $c_i$ refers to the center and $\sigma_i$ refers to the half-length of the interval of the $i^{th}$ variable, the function to optimize becomes:
\\
\begin{equation}
    J(c_1,..,c_n,\sigma_1,..,\sigma_n) =  \frac{1}{\prod_{i=1}^{n} 2\sigma_i} \times \int_{c_1 - \sigma_1}^{c_1+\sigma_1} .. \int_{c_n - \sigma_n}^{c_n+\sigma_n} f(x_1,..,x_n) dx_1..dx_n
\end{equation} 
\\
The problem with this objective function is that by optimizing it we indeed find a domain where our function f takes a large mean but no constraints are applied which means that the domain found could be very small, so we will add a gain that increases when the length of the intervals increase to give a bigger importance to large intervals ie: 

\begin{equation} 
 J(c_1,..,c_n,\sigma_1,..,\sigma_n) = \frac{1}{\prod_{i=1}^{n} 2 \sigma_i} \times \int_{c_1 - \sigma_1}^{c_1+\sigma_1} .. \int_{c_n - \sigma_n}^{c_n+\sigma_n} f(x_1,..,x_n) dx_1..dx_n  + \lambda \sum_{i=1}^{n} \sigma_i^2
\end{equation}

\subsection{Optimizing the objective function}

Now that we have defined our objective function we have to optimize it, the first approach would be to directly use the gradient ascent algorithm by computing the gradient of J at each iteration ie : 
\begin{equation*}
    c_i \leftarrow c_i + \alpha \times \frac{\partial J}{\partial c_i} 
    \quad
    \sigma_i \leftarrow \sigma_i + \alpha \times \frac{\partial J}{\partial \sigma_i} 
\end{equation*}
The problem is that in the vast majority of cases we can't compute the expression of that integral because of the complexity of the function and so we can't compute the gradient.
\newline So here the idea is to approximate at each iteration of the optimization process our function f by another function g, any function g that satisfies two proprieties can be used, first, it should be able to approximate f on the current domain found in the precedent iteration, second, we should be able to compute the expression of its integral.
\newline The family of the polynomial functions satisfies those two properties, for this paper we choose a polynomial function of degree 2 ie: 
\begin{equation}
    g(x_1,...,x_n) = \sum_{i=1}^{n} a_i x_i + \sum_{i=1}^{n} \sum_{j=i}^{n} b_{ij} x_ix_j  +d
\end{equation}
At each iteration, the coefficients $a_i, b_{ij}, d$ are found by:
\begin{enumerate}[label=(\roman*)]
 \item Sampling uniformly from the current domain K points (K has to be big enough ) $(x_{11},x_{21},..x_{n1},y_1) .. (x_{1K},x_{2K},..x_{nK},y_K)$ where $y_{i} = f(x_{1i},x_{2i},..x_{ni})$.
 \\
 \item Compute the variables $x_i x_j$ for each points.
 \\
 \item Perform a linear regression that minimize the mean squared error loss where the inputs are the points $(x_{1i},..,x_{ni},x_{1i}x_{1i},x_{1i}x_{2i},..,x_{ni}x_{ni})$ and the output is $y_i$.
\end{enumerate} 

At each iteration t of the optimization process we can approximate our objective function J by another objective function $L_t$ such as: 
\\
 \begin{equation}
    L_t(c_1,..,c_n,\sigma_1,..,\sigma_n) = 
    \frac{1}{\prod_{i=1}^{n} 2 \sigma_i} \times \int_{c_1 - \sigma_1}^{c_1+\sigma_1} .. \int_{c_n - \sigma_n}^{c_n+\sigma_n} g_t(x_1,..,x_n) dx_1..dx_n  + \lambda  \sum_{i=1}^{n} \sigma_i^2
\end{equation}
    g is integrable since it is a polynomial function, the expression of the integral of each term of g is given by: 
\\
\begin{itemize}[label=$\bullet$]

 \item $\int_{c_1 - \sigma_1}^{c_1+\sigma_1} .. \int_{c_n - \sigma_n}^{c_n+\sigma_n} d \quad dx_1..dx_n$ = d $\prod_{k=1}^{n} 2\sigma_k$
 
  \item $\int_{c_1 - \sigma_1}^{c_1+\sigma_1} .. \int_{c_n - \sigma_n}^{c_n+\sigma_n} a_i x_i  dx_1..dx_n$ = $a_i c_i \prod_{k=1}^{n} 2\sigma_k$
  
  \item $\int_{c_1 - \sigma_1}^{c_1+\sigma_1} .. \int_{c_n - \sigma_n}^{c_n+\sigma_n} b_{ij} x_i x_j  dx_1..dx_n$ = $b_{ij} c_i c_j \prod_{k=1}^{n} 2\sigma_k  \quad \text{if} \: i \ne j$
  \item $\int_{c_1 - \sigma_1}^{c_1+\sigma_1} .. \int_{c_n - \sigma_n}^{c_n+\sigma_n} b_{ij} x_i x_j  dx_1..dx_n$ = $b_{ij} (c_i^2 + \frac{\sigma_i^2}{3})  \prod_{k=1}^{n} 2\sigma_k \quad \text{if} \: i = j$
\end{itemize}

We then apply the gradient ascent algorithm ( or any other variant ) at each iteration  ie: 
\begin{equation*}
    c_i \leftarrow c_i + \alpha \times \frac{\partial L_t}{\partial c_i} 
    \quad
    \sigma_i \leftarrow \sigma_i + \alpha \times \frac{\partial L_t}{\partial \sigma_i} 
\end{equation*}

\subsection{Constraints}
What would be interesting to do is to not find the domain D from scratch but to set the values ( or the intervals ) of some variables and the algorithm will try to complete the others.
Formally if we want to set the value of the variable $x_i$ to for example 1, then the function h defined by: 
\begin{equation}
    h(x_1,..,x_{i-1},x_{i+1},..,x_m) = f(x_1,..,x_{i-1},1,x_{i+1},..,x_m)
\end{equation}
Will be used instead of f in the optimization process to complete the domain.
If now we want to set an interval for the variable $x_i$ which has 0 as center and 1 as half-length than the function h used instead of f becomes:
\begin{equation}
  h(x_1,..,x_{i-1},x_{i+1},..,x_m) = \int_{-1}^{1} f(x_1,..,x_{i-1},x_i,x_{i+1},..,x_m) dx_i
\end{equation}
In other words it would be equivalent to say that the values of  $c_i$ and $\sigma_i$ would be fixed to 0 and 1 respectively and won't be part of the optimization process.

\section{Use case in machine learning}
We can apply the iterative domain optimization algorithm to machine learning models, for example let's say that we have a model that performs a binary classification by taking $(x_1,..,x_n)$ features and output a probability q $\in$ [0,1], the goal is to find a domain $[a_1,b_1]\times .. \times [a_n,b_n]$ where the values of q are large or small.
\\
\newline Let's first suppose than we only deal with numerical features, problems that can occurs is that we find a domain that haven't been seen by the model during the training step where the behavior is completely random, to avoid that issue we suppose that the data used in the training step was normalized with the Standard Score so that the features are centered in 0 and then penalize intervals which center deviate from 0, the new objective function becomes:
\\
\begin{equation}
 J(c_1,..,c_n,\sigma_1,..,\sigma_n) = \frac{1}{\prod_{i=1}^{n} 2 \sigma_i} \times \int_{c_1 - \sigma_1}^{c_1+\sigma_1} .. \int_{c_n - \sigma_n}^{c_n+\sigma_n} f(x_1,..,x_n) dx_1,..,dx_n  \; + \lambda \sum_{i=1}^{n} \sigma_i^2 - \beta \sum_{i=1}^{n} c_i^2
\end{equation}
\\
Let's now suppose that both numerical and categorical features are used, generally when we deal with categorical features we perform some prepossessing before feeding the data to the model, we will suppose here that the one-hot-encoding is used to those features meaning that each categorical feature $x_i$ is split into $K_i$ new features representing the modalities.
\\
\newline We start enumerating the numerical features from 1 to n and the categorical features from n+1 to m, each categorical feature $x_i$ that had $K_i$ modalities will be split into $K_i$ new features
$x_i^1 .. x_i^{K_i}$.
\\
\newline Like before our algorithm will try to find intervals for each features except that some additional constraints must be satisfied, first, each new feature representing a modality can only takes the values 0 or 1 due to the encoding used, second, the sum of the values of all features representing the modalities of a categorical feature must be equal to 1.
\\
\newline Since the algorithm find intervals for each feature, we can force each interval to be centred in 0 or 1 and have a very small length by adding two penalities terms, also the sum of the centers of the intervals of the features representing the modalities of a categorical feature can be close to 1  with a third penality term. 
Our objective function becomes:

\begin{equation}
\begin{split}
&\hspace*{-40pt}J(c_1,..,c_n,c_{n+1}^1,..,c_{m}^{K_{m}},\sigma_1,..,\sigma_n,\sigma_{n+1}^1,..,\sigma_m^{K_{m}}) =\\
&\hspace*{-40pt}\dfrac{1}{(\prod\limits_{i=1}^{n} 2 \sigma_i)
(\prod\limits_{i=n+1}^{m} \prod\limits_{p=1}^{K_i} 2\sigma_i^{p})} \times \int_{c_1 - \sigma_1}^{c_1+\sigma_1} .. \int_{c_m^{K_m} - \sigma_m^{K_m}}^{c_m^{K_m}+\sigma_m^{K_m}} f(x_1,..,x_m^{K_m}) dx_1..dx_m^{K_m} \\ 
&\hspace*{-50pt}+ \lambda \sum\limits_{i=1}^{n} \sigma_i^2 - \beta  \sum\limits_{i=1}^{n} c_i^2 \\
&\hspace*{-50pt}- \mu \sum\limits_{i=n+1}^m \sum\limits_{p=1}^{K_i} (c_i^{p}(1-c_i^{p}))^2 
-\omega \sum\limits_{i=n+1}^m \sum\limits_{p=1}^{K_i} (\sigma_i^{p})^2  
- \gamma \sum\limits_{i=n+1}^{m} ( \sum\limits_{p=1}^{K_i} c_i^{p} - 1)^2
\end{split}
\end{equation}

\begin{enumerate}[label=(\roman*)]
\item The first new penality term penalize when the center of the intervals of the features representing the modalities deviate from both 1 and 0.
\\
\item The second new penality term penalize the large intervals of the features representing the modalities.
\\
\item The last new penality term penalize when the sum of the center of the intervals of the features representing the modalities of a categorical feature deviate from 1.
\end{enumerate}
 
 Note that those new features are not normalized ( else it would be difficult since we have to keep track of the normalized values of 0 and 1 for each features ), also the new tuning parameters have to be large enough due to the importance of the new constraints.

\section{Results}
For the experiments we will use the titanic dataset and make some preprocessing before training two models : a neural network and a random forest.
We fill the missing values of the age feature, we drop unnecessery features and keep the age, parsh, sibsp, fare, sex, pclass and embarked features. After that we normalize the numericals features with the standard scores and we one hot encode the categorical features, we end up with 12 features.
We then divide our dataset into 80\% training and 20\% test.
the following table summarize the results on the test set. Note that since it's a classification task the models take their values in [0,1].
\\
\setlength{\tabcolsep}{0.5em}
\begin{center}
\begin{tabular}{|c|c|c|}
\hline 
Model & auc & accuracy \\
\hline 
neural network & 0.89 & 0.84\\
\hline 
random forest & 0.84 & 0.8 \\
\hline

\end{tabular}
\end{center}

We will see the four following cases:
\begin{enumerate}[label=(\roman*)]
 \item Optimizing a domain where the neural network takes large values.
 \item Optimizing a domain where the random forest takes large values. 
 \item Optimizing a domain where the neural network takes small values. 
 \item Optimizing a domain where the neural network takes large values with two constraints : the modality of the embarked feature is set to Q and the normlized value of the age feature  $\in [1,3]$.
\end{enumerate}
A last case will be explored with an estimation of the probability density function of the input features using the kernel density estimation ie: 
\begin{equation}
    \hat{p}_H(x) = \frac{1}{n} \sum_{i=1}^{n} K_H(x-x_i)
\end{equation}
where the kernel used is gaussian and a bandwidth of 0.2, this estimator is trained on the whole dataset and 
in this case the integral to optimize is:
\begin{equation}
    \frac{1}{Vol(D)} \times \int_D f(x_1,..,x_n) \hat{p}_H(x_1,..,x_n)
\end{equation}
\\
The values of the tunning parameters are summarized in the follow table :
\begin{center}
\begin{tabular}{c|c|c|c|c|c|c|c}
     \hline
     Cases & K & learning rate & $\lambda$ & $\beta$ & $\mu$ & $\omega$ & $\gamma$  
     \\
     \hline
     Case 1 & 50 & 0.07 & 0.1 & 0.03 & 30 & 1 & 10  \\ 
     \hline
     Case 2 & 50 & 0.07 & 0.001 & 0.005 & 3 &1 & 3  \\ 
     \hline  
     Case 3 & 50 & 0.07 & 0.1 & 0.03 & 30 & 1 & 10  \\ 
     \hline
     Case 4 & 50 & 0.07 & 0.1 & 0.01 & 10 & 1 & 10  \\ 
     \hline
     Case 5 & 80 & 0.2& 0.003 & 0.003 & 10 & 0.1 & 10  \\ 
     \hline     
\end{tabular}
\end{center}
Each case is launched for 300 iterations using the Adam optimizer,  surprisingly a polynomial function of degree 2 was enough to achieve good results, they are summarized in figure 1.

\begin{figure}
    \centering
     \includegraphics[height=8cm,width=14cm]{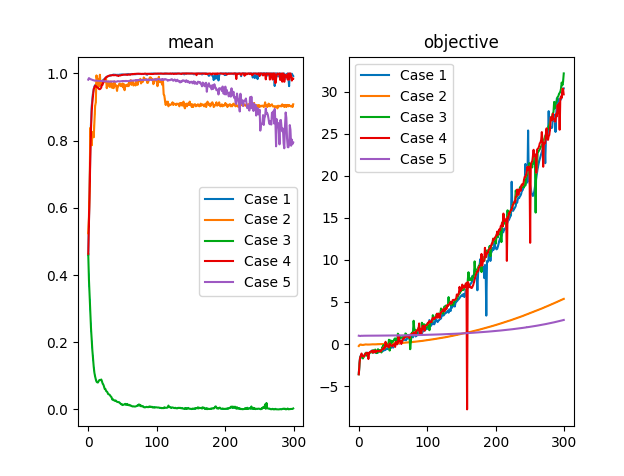}
    \caption{Evolution of the mean of the model and the objective function over the iterations for each case. }
    \label{fig:my_label}
\end{figure}
For the case with the neural network without the kernel density estimation (case 1, 3 and 4) we can see that the mean is fastly improved until it reaches its maximum or minimum while the objective function keep increasing, the large values of the objective function are due to the large domains found because of the large values of the $\sigma_i$. At the end of the optimization process all the constraints were satisfied for both numerical and categorical features.
\\
Case 2 with the random forest needed a small values for  $\lambda$ and $\beta$ to stabilize the optimization, the consequence of that is a small value for the objective function compared to cases 1, 3 and 4. Similarly to the previous cases the mean reach fastly its maximum while the objective function is still improving, the constraint are also satisfied at the end. 
\\
Case 5 with the kernel density estimator is particular, here we initialized the optimization process with a small domain centred on a row of the dataset because the density function applied on a random initialization will value 0 which can causes some difficulties. The values of $\beta$ and $\lambda$ are small for the same reason as case 2 and at the end all the constraints were also satisfied.
\\
The optimization method give satisfying results even for a complex function to optimize like case 5.
\section{Conclusion and future work}
The main contribution of this paper was to show the efficiency of the iterative domain optimization algorithm to find domains that maximize or minimize the values of a given function. We explored a use case of this algorithm in machine learning to 'interpret' models, However, a good estimator of the probability density function of the inputs features is needed to end up with coherent domains. We hope that this contribution will help in the optimization area. For future works we aim to improve the optimization algorithm and apply it in others fields.

\end{document}